\documentclass[12pt]{amsart}
\usepackage{amssymb}

\theoremstyle{plain}
\newtheorem*{theorem}{Theorem}

\newtheorem*{definition}{Definition}

\newtheorem*{remark}{Remark}

\newcounter{passvalue}

\newcommand{\reff}[1]{(\ref{#1})}
\newcommand{\ds}{\displaystyle}
\newcommand{\cstar}{\ensuremath{\text{C}^{*}}\nobreakdash-\hspace{0 pt}}
\newcommand{\BB}{\mathcal{B}}
\newcommand{\HH}{\mathcal{H}}
\newcommand{\LL}{\mathcal{L}}
\newcommand{\MM}{\mathcal{M}}
\newcommand{\alg}{\operatorname{Alg}}
\newcommand{\mx}{\mathcal{M}_x}

\begin{document}
\title[Invariant Linear Manifolds]{Invariant Linear Manifolds\\
        for CSL-algebras and Nest Algebras}
\author{Alan Hopenwasser}
\address{Department of Mathematics\\
        University of Alabama\\
        Tuscaloosa, AL 35487}
\email{ahopenwa@euler.math.ua.edu}
 \thanks{The author would like to thank Ken Davidson
 for drawing to his attention the references regarding
 operator ranges.}

 \date{May 18, 1999}

 \begin{abstract}
Every invariant linear manifold for a CSL-algebra,
$\operatorname{Alg} \mathcal{L}$, is a closed subspace
if, and only if, each non-zero projection
in $\mathcal{L}$ is generated by finitely many atoms 
associated with the projection lattice.  When $\mathcal{L}$
is a nest, this condition is equivalent to the condition that
every non-zero projection in $\mathcal{L}$ has an immediate
predecessor ($\mathcal{L}^{\perp}$ is well ordered). 
 The invariant linear manifolds of a nest algebra are
totally ordered by inclusion if, and only if, every non-zero
projection in the nest has an immediate predecessor.
 \end{abstract}

\maketitle

 Kadison's transitivity theorem \cite{rvk57} implies that a \cstar
 algebra acting on a Hilbert space which has no non-trivial closed
 invariant subspaces must also have no non-trivial invariant linear
 manifolds.  This note investigates the analogous situation for
 CSL-algebras and, in particular, for nest algebras.  We identify
 exactly the family of CSL-algebras for which every invariant linear
 manifold is, in fact, a closed subspace (and hence an element of the
 lattice of invariant subspaces for the algebra).  When specialized to
 nest algebras, this family reduces to those nest algebras for which
 each non-zero element of the nest has an immediate predecessor.
 Note that a nest, $\LL$, satisfies this condition if, and only if,
 the nest, $\LL^{\perp}$, consisting of orthogonal complements of
 elements of $\LL$ is well ordered.  For
 any other nest algebra, not only are there invariant linear manifolds
 which are not closed, but the family of invariant linear manifolds is
 not totally ordered by inclusion.  Thus, most nest algebras are
 ``nest algebras'' in a topological sense but not an algebraic sense.

Some of the results in this paper were obtained by Foia\c{s} in
1971/72 in his study of invariant operator ranges \cite{cp71, cp72} or
by Ong in his description of all the invariant operator ranges for a
nest algebra \cite{sco80}. Davidson extended this description to
CSL-algebras in \cite{krd82}.  Davidson's book on nest algebras
\cite{krd88bk} is a convenient reference for most of this work on
operator ranges.  All the same, we give complete and independent
proofs for those of our results on arbitrary linear manifolds which
follow from prior work on operator ranges; these proofs are shorter
and more elementary than the arguments via operator ranges.  (The
price to be paid is that these results are somewhat weaker than the
operator range results.  The converses, on the other hand, are
correspondingly stronger.)

To fix notation, let $\LL$ be a commutative lattice of orthogonal
projections acting on a Hilbert space $\HH$.  It is assumed that $\LL$
is complete and that it contains $0$ and $I$.  The algebra of all
bounded linear operators on $\HH$ which leave invariant each
projection in $\LL$ is denoted by $\alg \LL$.  It is convenient to use
the same symbol to denote both a subspace and the orthogonal
projection whose range is the subspace.  The term ``subspace'' will
always mean ``closed linear subspace''; the term `linear manifold' is
used when closure is not assumed.  Thus $\LL$ will be viewed
simultaneously as a 
complete lattice of closed subspaces of $\HH$ and as a
lattice of mutually commuting 
projections in $\BB(\HH)$ which is closed
in the strong operator topology (and which contains $0$ and $I$).

The primary tool used in this note is the necessary and sufficient
condition on two vectors $x$ and $y$ in $\HH$ for the existence
of an operator $T$ in $\alg \LL$ such that $Tx=y$: 
with the understanding
that the fraction $0/0$ is to be interpreted as $0$, the existence
of an operator $T$ carrying $x$ to $y$ is equivalent to
\[
\sup_{E\in \LL} \frac {\Vert E^{\perp}y \Vert}
  {\Vert  E^{\perp}x \Vert} < \infty.
\]
This was first proved by Lance \cite{ecl69} for nest
algebras and then extended to CSL-algebras in \cite{ah80}.
 
Recall that an \textit{atom}, $A$, associated with a subspace lattice
$\LL$ is a minimal non-zero interval from $\LL$; i.e.,
$A$ has the form $P-Q$ where $P,Q \in \LL$, $Q < P$, and
$AF$ is either $0$ or $A$, for all $F \in \LL$.  If
$A_1, A_2, \dots $ is a sequence of atoms from $\LL$ (either
a finite or an infinite sequence)  let $E(A_1, A_2, \dots)$
denote the smallest projection in $\LL$ which contains each of
the atoms $A_1, A_2, \dots$.  In other words,
\[
E(A_1, A_2, \dots ) = \bigwedge_{F\in \LL}\{F \mid A_n \subseteq F,
\text{ for all $n$}\}.
\]
In the following definition, and throughout this paper, when we 
say that a projection in $\LL$ is generated by a set of atoms, we 
mean that it is the smallest projection in $\LL$ which contains each
of the atoms of the set.

\begin{definition}
A commutative subspace lattice, $\LL$, is {\upshape hyperatomic} if
every non-zero projection in $\LL$ is generated by finitely many atoms.
\end{definition}

\begin{remark}
If $P = E(A_1, \dots ,A_n)$, then we may as well assume that
the atoms $A_1, \dots A_n$ are independent in the sense that
$A_iE(A_j) = 0$ whenever $i \ne j$.  (Just delete some atoms
from the list, if necessary.)
\end{remark}

If $P$ is a non-zero projection in $\LL$, let $P_{-}$ denote the
projection $\bigvee \{F\in \LL \mid F \ngeqq P \}$.  If $x$ and $y$
are vectors in $\HH$, then there is a rank-one operator
$T \in \alg \LL$ such that $Tx=y$ if, and only if, 
$x \in P_{-}^{\perp}$ and $y \in P$ for some $P \in \LL$.
For nests, this was proven by Ringrose \cite{jr65}; the extension
to commutative lattices is due to Longstaff \cite{wl75}.
Note that if $A$ is an atom from $L$ and if $P =E(A)$, then
$A \leq P_{-}^{\perp}$.  Consequently, if $x \in A$ and
$y \in P$, then there is a rank-one operator in $\alg \LL$
such that $Tx=y$. 

In the theorem which follows, \textit{operator range} refers
to the range of a bounded linear operator acting on $\HH$.  

\begin{theorem}
Let $\LL$ be a commutative subspace lattice acting on a separable
Hilbert space.
The following conditions are equivalent:
\begin{enumerate}
\item \label{hyper} $\LL$ is hyperatomic.
\item \label{ascend} Every ascending sequence, 
$F_1 \leq F_2 \leq \dots$, of projections
in $\LL$ is eventually constant. 
\item \label{oprange} Every invariant operator
range for $\alg \LL$ is a closed subspace (and therefore
an element of $\LL$).
\item \label{manifold} Every invariant linear manifold
for $\alg \LL$ is a closed subspace (and therefore
an element of $\LL$).
\setcounter{passvalue}{\value{enumi}}
\end{enumerate}
\noindent If $\LL$ is a nest, then the following additional conditions
are  equivalent to each of the conditions above:
\begin{enumerate}
\setcounter{enumi}{\value{passvalue}}
\item \label{pred} Every non-zero projection in $\LL$ has
an immediate predecessor.
\item \label{oprangenest} The invariant operator ranges for
$\alg \LL$ are totally ordered by inclusion.
\item \label{manifoldnest} The invariant linear manifolds for
$\alg \LL$ are totally ordered by inclusion.
\end{enumerate}
\end{theorem}
\begin{proof}
\reff{hyper}$\Rightarrow$\reff{ascend}: Let $F_1 \leq F_2 \leq \dots$
be an ascending sequence of projections in $\LL$.  Let $P = \bigvee
F_j$.  Since $\LL$ is hyperatomic, there are finitely many atoms,
$A_1, \dots, A_n$, so that $P = E(A_1, \dots A_n)$.  Each of these
atoms must be a subprojection of some $F_j$; since there are only
finitely many atoms in this list, there is $j$ such that $A_i \leq
F_j$, for all $i = 1, \dots n$.  But then $P = F_j$; hence $P = F_k$,
for all $k \geq j$.

\reff{ascend}$\Rightarrow$\reff{hyper}: Assume that $\LL$ is not
hyperatomic.  Let $P \in \LL$ be a projection which is not generated
by finitely many atoms.  Let $Q$ be the smallest projection in $\LL$
which contains all atoms which are 
subprojections of $P$.  First, suppose
that $Q < P$.  If $F$ is any element of $\LL$ such
that $Q \leq F < P$, then there is a projection $G \in \LL$ such that
$F < G < P$ (since $P-F$ is not an atom from $\LL$).  A routine induction
argument now yields an ascending sequence of projections which is not
eventually constant.  Next, suppose that $Q = P$.  In other words,
there is a sequence $A_1, A_2, \dots$ (necessarily infinite) of atoms
such that $P = E(A_1, A_2, \dots)$.  For each $n$, let 
$F_n = E(A_1, \dots, A_n)$.  Then $F_1, F_2, \dots$ is an ascending
sequence of projections in $\LL$, $F_n < P$ for all $n$, and
$\bigvee F_n = P$.  Thus, this sequence is not eventually constant.

\reff{manifold}$\Rightarrow$\reff{oprange} is immediate and
\reff{oprange}$\Rightarrow$\reff{ascend} follows from
\cite[Theorem 15.29]{krd88bk}.  However, we give here an
elementary proof of \reff{manifold}$\Rightarrow$\reff{ascend} which does
not require any information about operator ranges.
Suppose \reff{ascend} is false.  Let $0 < F_1 < F_2 < \dots$ be a
strictly increasing sequence of projections in $\LL$.  For each $n$,
let $x_n \in F_n - F_{n-1}$ be a non-zero vector chosen so that
$\ds \sum_{n=1}^{\infty}n^2 \Vert x_n \Vert^2 < \infty$.  (For
example, let $x_n$ be any vector in $F_n - F_{n-1}$ for which
$\Vert x_n \Vert = 1/n^2$.)

For each $n$, let $y_n = nx_n$.  By the choice of the $x_n$, the
sequence, $y_n$, is a square summable sequence of mutually 
orthogonal vectors.  Let $\ds x = \sum_{n=1}^{\infty}x_n$ and
$\ds y = \sum_{n=1}^{\infty}y_n$.  (Both sums converge in $\HH$.)
Let $\mx = \{Tx \mid T \in \alg \LL \}$, an invariant linear 
manifold for $\alg \LL$.  For each $n$,
\begin{align*}
F_n^{\perp}x &= \sum_{k=n+1}^{\infty}x_k \quad \text{and}\\
F_n^{\perp}y &= \sum_{k=n+1}^{\infty}y_k = \sum_{k=n+1}^{\infty}
  kx_k.
\end{align*}
So,
\begin{align*}
\Vert F_n^{\perp}y \Vert^2 
   &= \sum_{k=n+1}^{\infty}k^2 \Vert x_k \Vert ^2
        \geq \sum_{k=n+1}^{\infty} (n+1)^2 \Vert x_k \Vert^2 \\
   &= (n+1)^2 \sum_{k=n+1}^{\infty} \Vert x_k \Vert^2
        = (n+1)^2 \Vert F_n^{\perp}x \Vert^2
\end{align*}
Thus,
\[
\frac {\Vert F_n^{\perp}y \Vert}{\Vert F_n^{\perp}x \Vert}
\geq n+1
\]
and hence,
\[
\sup_{F \in \LL}\frac {\Vert F^{\perp}y \Vert}{\Vert F^{\perp}x \Vert}
= \infty.
\]
This shows that $y \notin \mx$.

On the other hand, $(F_n - F_{n-1})x = x_n \in \mx$, for all $n$,
whence $y_n = nx_n \in \mx$ and 
$\ds \sum_{n=1}^K y_n \in \mx$, for all $K$.  Since
$\ds y = \lim_{K \to \infty} \sum_{n=1}^K y_n$, this shows that
$y \in \overline{\mx}$.  Thus, $\mx$ is an invariant linear manifold 
which is not closed.

\reff{hyper}$\Rightarrow$\reff{manifold}: First, we show that if
$x \in \HH$, then $\mx$ is closed; i.e., every singly generated
invariant linear manifold is closed.  
        Let $x \in \HH$ and $P = \overline{\mx}$.  If $x=0$, then
$P=0$; so assume $x \ne 0$ and, hence, $P \ne 0$.  Since $\LL$
is hyperatomic, there exist atoms $A_1, \dots, A_n$ so that 
$P = E(A_1, \dots, A_n)$.  Without loss of generality, we may assume
that $A_1, \dots , A_n$ are independent; i.e., that $A_iE(A_j) = 0$
whenever $i \ne j$.  Equivalently, $A_iTA_j = 0$ for all
$T \in \alg \LL$, when $i \ne j$.  
It follows that $A_jx \ne 0$, for all $j$.

If $y \in E(A_j)$, then, by the comments preceding the statement
of the theorem, there is $T \in \alg \LL$ such that 
$y = TA_jx$.  Thus $E(A_j) \subseteq \mx$, for all $j$.

 If $y \in P$ is arbitrary, then, since 
$P = E(A_1) \vee \dots \vee E(A_n)$, there exist $y_i \in E(A_i)$
such that $y= y_1 + \dots y_n$.  Since each $y_i \in \mx$, we
have $y \in \mx$.  This shows that 
$P =\overline{\mx} \subseteq \mx$,
so $\mx$ is closed.

Before turning to general invariant linear manifolds, we need an
observation:  if $P_1 = \MM_{x_1}$ and $P_2 = \MM_{x_2}$, then
there is a vector $x$ such that $\mx = P_1 \vee P_2$.  Indeed,
choose $x = x_1 + P_1^{\perp}x_2$.  Since $x_1 = P_1x$, we have
$P_1 \subseteq \mx$.  Let $y \in P_2P_1^{\perp}$.  Since
$P_2 = \MM_{x_2}$, there is $T \in \alg \LL$ such that 
$y = Tx_2$.  Since $x_2 = P_1^{\perp}x_2 + P_1 x_2 $,
\[
y = Tx_2 = TP_1^{\perp}x_2 + TP_1x_2 = TP_1^{\perp}x_2 
 + P_1TP_1x_2.
\]
But $P_1^{\perp}y = y$, so
\[
y = P_1^{\perp}TP_1^{\perp}x_2 = P_1^{\perp}TP_1^{\perp}
(x_1 + P_1^{\perp}x_2) = P_1^{\perp}TP_1^{\perp}x.
\]
Thus, $P_2P_1^{\perp} \subseteq \mx$.  Since $P_1 \vee P_2
 = P_1 + P_2P_1^{\perp}$, we have $P_1 \vee P_2 \subseteq \mx
 \subseteq P_1 \vee P_2$; i.e., $\mx = P_1 \vee P_2$.

Finally, let $\MM$ be an arbitrary invariant linear manifold
for $\alg \LL$.  Let $P = \overline{\MM}$.  So
$P = \bigvee \{ \mx \mid x \in \MM \}$.  Since $\HH$
is separable, we can write $P$ as the join of countably many
subspaces of the form $\mx$; i.e., there is a sequence,
$x_1, x_2, \dots$, of vectors in $\MM$ so that $\ds P =
\bigvee_{j=1}^{\infty} P_j$, where $P_j = \MM_{x_j}$, all $j$.
The observation above shows that $P_1 \vee P_2$ and, indeed,
 any finite join $P_1 \vee \dots \vee P_n$ can be written in
the form $\mx$ for some $x \in \MM$.  So we may assume that
$P_1 \leq P_2 \leq \dots$.  Since $\LL$ is hyperfinite and we have
already shown that \reff{hyper}$\Rightarrow$\reff{ascend}, this sequence
is eventually constant.  But this shows that $P = \MM_{x_j}$, for some
$x_j \in \MM$.  Thus $P = \MM$ and $\MM$ is closed.

For the rest of the proof we assume that $\LL$ is a  nest.

The equivalence, \reff{hyper}$\Leftrightarrow$\reff{pred}, is
trivial; in a nest a projection is generated by finitely many
atoms if, and only if, it is generated by a single atom.
\reff{hyper}$\Rightarrow$\reff{manifoldnest} follows immediately
from \reff{hyper}$\Rightarrow$\reff{manifold} and 
\reff{manifoldnest}$\Rightarrow$\reff{oprangenest} is trivial.

 \reff{oprangenest}$\Rightarrow$\reff{ascend}:  Assume that
\reff{ascend} if false; let $0 = F_0 < F_1 < F_2 < \dots$
be a strictly increasing sequence of projections in $\LL$.
Suppose that $\lambda = (\lambda_n)$ is a decreasing sequence
of positive real numbers.  Let $\ds D_{\lambda} =
\sum_{n=1}^{\infty} \lambda_n (F_n - F_{n-1}) $ (the sum
converges in the strong operator topology).  By a result
of Ong \cite{sco80}, the range of $D_{\lambda}$
is an invariant operator range for $\alg \LL$.

If $x$ is a vector in $\ds \bigvee_{n=1}^{\infty}F_n$, let
$x_n = (F_n - F_{n-1})x$, for each $n$.  Necessarily,
$\ds \sum_{n=1}^{\infty}\Vert x_n \Vert^2 < \infty$.
Note that $x$ is in the range of $D_{\lambda}$ if, and only
if, $\ds \sum_{n=1}^{\infty} \frac 1{{\lambda_n}^2}
\|x_n\|^2 < \infty$.

In order to exhibit two operator ranges which are not related by
inclusion, choose two decreasing seqences, $\lambda$ and $\mu$,
of positive real numbers in such a way that
$ {\mu_n}/{\lambda_n} \geq n$ whenever  $n$ is even and
$ {\lambda_n}/{\mu_n} \geq n$ whenever $n$ is odd.  Next,
choose vectors $x_n$ and $y_n$ in $F_n - F_{n-1}$ so that
$\Vert x_n \Vert =  {\mu_n}/n$ and
$\Vert y_n \Vert =  {\lambda_n}/n$, for all $n$.
The four sequences, $(x_n)$, $\ds \left( \frac 
 1{\mu_n}x_n \right)$, $(y_n)$, and
$\ds \left( \frac 1{\lambda_n}y_n \right)$ are all square
summable.  Consequently, $\ds x = \sum_{n=1}^{\infty}
x_n$ and $\ds y = \sum_{n=1}^{\infty}y_n$ are convergent
sums; $x$ is an element of the range of $D_{\mu}$; and
$y$ is an element of the range of $D_{\lambda}$.

Now consider the sequences $\ds \left( \frac 1{\lambda_n}
x_n \right)$ and $\ds \left( \frac 1{\mu_n}y_n  \right)$.
For $n$ even,
\[
\left\| \frac {x_n}{\lambda_n} \right\| = 
\left(\frac {\mu_n}{\lambda_n} \right)
\left( \frac 1{\mu_n} \right) \Vert x_n \Vert \geq 1;
\]
while, for $n$ odd,
\[
\left\| \frac {y_n}{\mu_n} \right\| = 
\left(\frac {\lambda_n}{\mu_n} \right)
\left( \frac 1{\lambda_n} \right) \Vert y_n \Vert \geq 1.
\]
This shows that $x$ is not in the range of $D_{\lambda}$
and $y$ is not in the range of $D_{\mu}$;  thus $D_{\lambda}$
and $D_{\mu}$ are not ordered by inclusion.

\reff{manifoldnest}$\Rightarrow$\reff{ascend}: This is evident
from the preceeding argument and the trivial implication
\reff{manifoldnest}$\Rightarrow$\reff{oprangenest}.
However, here is an alternative proof 
which avoids the use of operator ranges and which exhibits
two singly generated invariant linear manifolds which are
unrelated by inclusion.

Let $\ds e_n = 1 + 2 + \dots + n$, for all
$n$.  Since $e_n = e_{n-1} + n$, 
we have, for all $k$,  
\begin{align*}
 \sum_{n = k+1}^{\infty} \frac 1 {2^{e_n}} 
   &< \sum_{n=0}^{\infty} \frac 1 {2^{e_{k+1}+n}}
     =\sum_{n=0}^{\infty} \frac 1 {2^{e_k +k+1+n}}   \\
   &= \frac 1 {2^{e_k}} \sum_{n=0}^{\infty} \frac 1 {2^{k+1+n}}
      = \frac 1 { 2^{e_k} 2^{k} }.
\end{align*}
Define two sequences:
\begin{align*}
 a_n &=
\begin{cases}
  0,                 &\text{if $n$ is even;} \\
   1/{2^{e_n}}, &\text{if $n$ is odd;}
\end{cases}
\\
 b_n &=
\begin{cases}
   1/{2^{e_n}},  &\text{if $n$ is even;} \\
  0,                  &\text{if $n$ is odd.} 
\end{cases}
\end{align*}
In other words,
\begin{align*}
a &= \left(\frac 1 {2^{e_1}}, 0, \frac 1 {2^{e_3}},
   0, \frac 1 {2^{e_5}}, 0, \dots  \right), \\
b &= \left(0, \frac 1 {2^{e_2}}, 0, \frac 1 {2^{e_4}},
    0, \frac 1 {2^{e_6}}, \dots \right).
\end{align*}

Now assume that $k$ is an odd integer.  Then
\[
\sum_{n=k}^{\infty} b_n = \sum_{n=k+1}^{\infty} b_n
 < \sum_{n=k+1}^{\infty} \frac 1 {2^{e_n}}
 < \frac 1 {2^{e_k}2^{k}}.
\]
Therefore,
\[
\left( \sum_{n=k}^{\infty} b_n \right)^{-1} > 2^{k}2^{e_k}.
\]
Since we also have
\[
\sum_{n=k}^{\infty} a_n > \frac 1 {2^{e_k}},
\]
we obtain
\[
 \left( \sum_{n=k}^{\infty} a_n \right)
  \left( \sum_{n=k}^{\infty} b_n \right)^{-1}
  > 2^{k}, \quad \text{for all odd $k$.}
\]
Similarly,
\[
 \left(\sum_{n=k}^{\infty}b_n \right)
   \left(\sum_{n=k}^{\infty}a_n  \right)^{-1} > 2^{k},
  \quad \text{for all even $k$.}
\]
Thus we have
\begin{align*}
 &\sup_k \left(\sum_{n=k}^{\infty}a_n \right)
   \left(\sum_{n=k}^{\infty}b_n  \right)^{-1} = \infty, 
  \quad \text{and}\\
 &\sup_k \left(\sum_{n=k}^{\infty}b_n \right)
   \left(\sum_{n=k}^{\infty}a_n  \right)^{-1} = \infty.
\end{align*}
Assume that $\LL$ fails to satisfy \reff{ascend} and
let $0 = F_0 < F_1 < F_2 < \dots$ be a strictly
increasing sequence of projections in $\alg \LL$.
For each $n \geq 1$, choose vectors $x_n$ and $y_n$
in $F_n - F_{n-1}$ so that
$\ds \Vert x_n \Vert^2 = a_n $ and  
$\ds \Vert y_n \Vert^2 = b_n$;
finally, let $\ds x = \sum_{n=1}^{\infty}x_n$
and $\ds y = \sum_{n=1}^{\infty}y_n$.
Since $\ds \Vert F_k^{\perp}x \Vert^2 = \sum_{n=k+1}^{\infty}a_n$
 and $\ds \Vert F_k^{\perp}y \Vert^2 = \sum_{n=k+1}^{\infty}b_n$,
it follows that
\[
\sup_k \frac {\Vert F_k^{\perp}x \Vert}
  {\Vert F_k^{\perp}y \Vert} = \infty \quad \text{and} \quad
\sup_k \frac {\Vert F_k^{\perp}y \Vert}
  {\Vert F_k^{\perp}x \Vert} = \infty.
\]
Thus
\[
\sup_{F \in \LL} \frac {\Vert F^{\perp}x \Vert}
  {\Vert F^{\perp}y \Vert} = \infty \quad \text{and} \quad
\sup_{F \in \LL} \frac {\Vert F^{\perp}y \Vert}
  {\Vert F^{\perp}x \Vert} = \infty.
\]

This shows that for all $T \in \alg \LL$, $Ty \ne x$ and
$Tx \ne y$.  The linear manifolds
$\mx$ and
$\MM_y$
are  invariant 
under $\alg \LL$ and, since $x \notin \MM_y$ and $y \notin \mx$,
we have $\MM_x \nsubseteq \MM_y$ and $\MM_y \nsubseteq \MM_x$.
\end{proof}

\providecommand{\bysame}{\leavevmode\hbox to3em{\hrulefill}\thinspace}

\end{document}